\documentclass[reqno]{amsart}

\usepackage[all]{xy}

\numberwithin{equation}{section}

\theoremstyle{plain}
\newtheorem{lemma}{Lemma}[section]

\newtheorem{propo}[lemma]{Proposition}

\newtheorem{con}[lemma]{Conjecture}

\theoremstyle{definition}

\theoremstyle{remark}

\newcommand{\bs}{{\rm Bs\thinspace}}
\newcommand{\p}{\mathbb{P}}

\newcommand{\ls}{{\mathcal L}}
\newcommand{\ms}{{\mathcal M}}

\newcommand{\ci}{{\mathcal I}}

\newcommand{\cre}{{\rm Cr\thinspace}}

\newcommand{\rmap}{\dashrightarrow}

\begin{document}

\title{A conjecture on special linear systems of $\p^3$}
\author{Antonio Laface}
\address{
Dipartimento di Matematica, Universit\`a degli Studi di Milano,
Via Saldini 50, \newline 20100 Milano, Italy }
\email{antonio.laface@unimi.it}

\author{Luca Ugaglia}
\address{
Dipartimento di Matematica, Universit\`a degli Studi di Milano,
Via Saldini 50, \newline 20100 Milano, Italy }
\email{luca.ugaglia@unimi.it}

\begin{abstract}
In this note we deal with linear systems of $\p^3$ through fat
points. We consider the behavior of these systems under a
cubo-cubic Cremona transformation that allows us to produce a
class of special systems which we conjecture to be the only ones.
\end{abstract}
\maketitle

\section{Introduction and statement of the problem}

In what follows we assume that the ground field is algebraically
closed of characteristic 0. Consider $r$ points in general
position on $\p^n$, to each one of them associate a natural number
$m_i$ called the {\em multiplicity} of the point. We will denote
by $\ls = \ls_n(d,m_1,\ldots, m_r)$ the linear system of
hypersurfaces of degree $d$ through the $r$ points with the given
multiplicities. Define the {\em virtual dimension} of the system
as $v(\ls) = \binom{d+n}{n}-\sum\binom{m_i+n-1}{n}-1$ and its {\em
expected dimension} by $e=\max\{v,-1\}$. Observe that
$e\leq\dim(\ls)$ and that the inequality may be strict if the
conditions imposed by the points are dependent. In this case we
say that the system is {\em special}. The problem of classifying
special systems has been completely solved in the case
$m_1=\ldots=m_r=2$ (see~\cite{alhi}) and it has been largely
studied for linear systems on the plane
(see~\cite{cm1,cm2,ev,mi}). The main conjecture on the structure
of special planar systems has been formulated in~\cite{ha,hi}. In
this note we report about some recent results in the case of
$\p^3$. In~\cite{lu1} we gave a counterexample to a conjecture
(see~\cite{ci}) about the structure of special linear systems of
$\p^3$. Starting from this idea in~\cite{lu2} we analyzed the
behavior of linear systems under a cubo-cubic Cremona
transformation of $\p^3$. This allowed us to construct a class of
special linear systems which we conjecture to be all the possible
ones.\\
Throughout this note we will denote by $X$ the blow up of $\p^3$
along the $r$ fixed points and by $E_i\cong\p^2$ the exceptional
divisors. From the Riemann-Roch formula on a smooth threefold, we
obtain that $v(\ls) = (\ls(\ls-K_X)(2\ls-K_X)+c_2(X) \ls)/12$. If
the linear system can be written as $\ls=F+\ms$, where $F$ is a
fixed divisor of $\ls$ and $\ms$ is the residual system, then the
above formula implies:
\begin{equation}\label{add}
v(\ls) = v(\ms) + v(F) + \frac{F\ms(\ls-K_X)}{2}.
\end{equation}

All the results described in this note can be found in \cite{lu2}.

\newpage

\section{Cubic Cremona transformations and linear systems}

It is possible to consider the behavior of linear systems under a
birational transformation of $\p^3$. In particular we need to
consider a transformation which sends linear systems through
points into systems through points. Consider the system
$\ls_3(3,2^4)$; by putting the four double points in the
fundamental ones, the associated rational map is $\cre:
(x_0:x_1:x_2:x_3) \rmap (x_0^{-1}:x_1^{-1}:x_2^{-1}:x_3^{-1}).$
This birational map induces the following action on a linear
system $\ls$:
\begin{eqnarray}\label{cre-a1}
\cre(\ls) & = & \ls_3(d+k,m_1+k,\ldots, m_4+k,m_5,\ldots,m_r),
\end{eqnarray}
where $k=2d-\sum_{i=1}^4 m_i$. By using this transformation, it is
easy to see that if $2d < m_1 + m_2 + m_3$, then the plane through
the first three points is a fixed component of the system. Observe
that $\dim \cre(\ls) = \dim \ls$ but in general the virtual
dimensions of the two systems may be different as stated in the
following:

\begin{propo}\label{vir-change}
Let $\ls$ be a linear system such that $2d \geq m_i + m_j + m_k$
for any choice of $\{i,j,k\}\subset\{1,2,3,4\}$ then
\begin{equation}\label{vc}
v(\cre(\ls)) - v(\ls) = \sum_{t_{ij} \geq
2}\binom{1+t_{ij}}{3}-\sum_{t_{ij}\leq -2}\binom{1-t_{ij}}{3},
\end{equation}
where $t_{ij}=m_i+m_j-d$.
\end{propo}

In particular this implies that $v(\cre(\ls))\geq v(\ls)$ if the
degree of $\cre(\ls)$ is smaller than the one of $\ls$. This means
that as far as $2d < m_1+m_2+m_3+m_4$ we can perform a Cremona
transformation decreasing the degree and the multiplicities of the
system. If at some step we get a system such that $2d < m_1 + m_2
+ m_3$, then we remove the plane from the base locus. After a
finite number of steps we get a system with $2d\geq
m_1+m_2+m_3+m_4$, which we say to be in {\em standard form}.

\section{Conjecture}

To each linear system $\ls$ we associate the $1$-cycle $
\Gamma(\ls):=\sum_{t_{ij}\geq 1} t_{ij}l_{ij}, $ where
$t_{ij}=m_i+m_j-d$ and $l_{ij}$ is the line through $p_i$ and
$p_j$. Observe that by definition
$H^0(\ls\otimes\ci_{\Gamma(\ls)}) = H^0(\ls)$, since each line
$l_{ij}\in\Gamma(\ls)$ is contained in $\bs(\ls)$ with
multiplicity at least $t_{ij}$.
\begin{propo}\label{vir=}
The relation between the Euler characteristic of the two sheafs is
given by:
\[
\chi(\ls\otimes\ci_{\Gamma(\ls)})=\chi(\ls) + \sum_{t_{ij}\geq
2}\binom{t_{ij}+1}{3}.
\]
\end{propo}
This implies that for each $\Gamma = \sum \alpha_{ij}l_{ij}$ with
$2\leq  \alpha_{ij}\leq t_{ij}$ we have $$ \dim\ls - v(\ls) \geq
\sum \binom{\alpha_{ij}+1}{3} - h^2(\ls\otimes\ci_{\Gamma}).$$  It
is possible to prove that if $\Gamma$ is just a multiple line then
$h^2(\ls\otimes\ci_{\Gamma})=0$ and the system is special. If the
system is in standard form, then the graph associated to
$\Gamma(\ls)$ can be only one of the following:
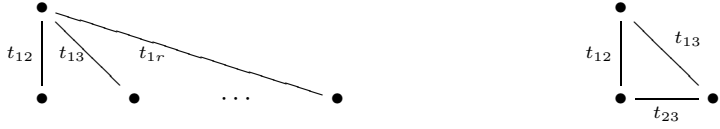
\[
\xymatrix{ \bullet \ar@{-}[d]_{t_{12}}^{\ t_{13}\hspace{7mm}t_{1r}} \ar@{-}[dr] \ar@{-}[drrr] & & & \\
\bullet & \bullet & \cdots & \bullet } \hspace{30mm}
\xymatrix{ \bullet \ar@{-}[d]_{t_{12}} \ar@{-}[dr]^{t_{13}} & & \\
\bullet \ar@{-}[r]_{t_{23}} & \bullet }
\]
The preceding discussion gives us a class of special systems in
standard form. We can construct another such class in the
following way. Let $Q=\ls_3(2,1^9)$ be a quadric through nine
points and suppose that $Q(\ls-Q)(\ls-K)<0$. By formula~\ref{add}
we have that $v(\ls) < v(\ls-Q)$ which implies that $\ls$ is
special. By assuming the Harbourne-Hirschowitz conjecture to be
true for planar systems with 10 points, it is possible to prove
that if $Q(\ls-Q)(\ls-K)<0$ then $Q$ is a fixed component of
$\ls$.

\begin{con}\label{conj}
A linear system $\ls$ in standard form is special if and only if
one of the following holds:
\begin{itemize}
\item[(i)] there exists a quadric $Q$ such that $Q(\ls-Q)(\ls-K)<0$;
\item[(ii)] at least one of the coefficients of $\Gamma(\ls)$ is
bigger then $1$.
\end{itemize}
\end{con}

Assuming that this conjecture and the Harbourne-Hirschowitz (for
systems through 10 points) are true we can remove all the quadrics
of step (i) from the base locus of $\ls$. Then the residual system
$\ls'$ is still in standard form and
$$\dim\ls = v(\ls')+\sum_{t'_{ij}\geq 2}\binom{t'_{ij}+1}{3}$$
assuming that $h^2(\ls'\otimes\ci_{\Gamma(\ls')})=0$.


We conclude this note with two propositions about {\em
homogeneous} linear systems, i.e. the systems $\ls$ for which
$m_1=\ldots=m_r=m$.
\begin{propo}
The system $\ls$ is empty for $d\leq 2m-1$ and $r\geq 8$.
\end{propo}
By assuming Conjecture~\ref{conj} and Harbourne-Hirschowiz
conjecture for linear systems on $\p^2$ with $10$ points, we can
also prove the following:
\begin{propo}
If $d\geq 2m$ the system $\ls$ is special if and only if $r=9$ and
$2m\leq d<[-1+\frac{3}{2}\sqrt{2m^2+2m}]$.
\end{propo}

Therefore if the system $\ls$ has more than $9$ fixed points (or
exactly $8$ points) then it is not special. If it has $9$ fixed
points, it is special if and only if $d$ satisfies the
inequalities of the preceding proposition. If $r\leq 7$ and $d\geq
2m$, the system can not be special. Finally, if $r\leq 7$ and
$d\leq 2m-1$, by applying a finite number of Cremona
transformations we reduce to a system in standard form.

\end{document}